\documentclass[12pt]{amsart}
\usepackage{latexsym, amsmath, amssymb}
\usepackage{epsfig}
\usepackage{amscd}
\usepackage{graphicx}
\usepackage[all]{xy}
\usepackage{pstricks,pst-node,cite}

\makeatletter
\newtheorem*{rep@theorem}{\rep@title}
\newcommand{\newreptheorem}[2]{%
\newenvironment{rep#1}[1]{%
 \def\rep@title{#2 \ref{##1}}%
 \begin{rep@theorem}}%
 {\end{rep@theorem}}}
\makeatother

\newtheorem{theorem}{Theorem}[section]
\newreptheorem{theorem}{Theorem}
\newtheorem{lem}[theorem]{Lemma}

\newtheorem{defi}[theorem]{Definition}
\theoremstyle{remark}

\psset{linewidth=.6pt, linecolor=blue, dash=3pt 3pt,doublesep=.05,dotsize=1pt 5}

\def\aa{{\mathcal A}}
\def\bb{{\mathcal B}}
\def\cc{{\mathcal C}}
\def\dd{{\mathcal D}}
\def\ee{{\mathcal E}}
\def\Aut{{\mathrm{Aut}}}
\def\Sym{{\mathrm{Sym}}}

\SpecialCoor

\ifx\blackandwhite\undefined
  \newrgbcolor{darkred}{.9 0 0}\newrgbcolor{emgreen}{0 .9 0}
  \newrgbcolor{pup}{.7  0 .9}\newrgbcolor{xxx}{.8 .8 .8} \else
\fi
\psset{linecolor=blue}
\psset{linewidth=.6pt, doublesep=.05,dotsize=1pt 5}
\newcommand{\middlearrow}{\lput{:U}{\pspicture[shift=0](0,0)(0,0)
\psline[arrows=->,arrowscale=1.5](2.2pt,0)(2.3pt,0)\endpspicture}}
\newcommand{\mda}{\lput{:U}{\pspicture[shift=0](0,0)(0,0)
\psline[linecolor=darkred,arrows=->,arrowscale=2.2](3.2pt,0)(3.4pt,0)\endpspicture}}
\newcommand{\mdb}{\lput{:U}{\pspicture[shift=0](0,0)(0,0)
\psline[linecolor=emgreen,arrows=->,arrowscale=1.7](3.2pt,0)(3.4pt,0)\endpspicture}}
\newcommand{\mdc}{\lput{:U}{\pspicture[shift=0](0,0)(0,0)
\psline[linecolor=pup,arrows=->,arrowscale=1.7](3.2pt,0)(3.4pt,0)\endpspicture}}

\begin{document}

\author{Dongseok Kim}
\address{Department of Mathematics \\Kyonggi University
\\ Suwon, 443-760 Korea}
\email{dongseok@kgu.ac.kr}

\title[A classification of transitive links]
{A classification of transitive links and periodic links}

\begin{abstract}
We generalized the periodic links to \emph{transitive} links in a $3$-manifold $M$.
We find a complete classification theorem of transitive links in a $3$-dimensional sphere $\mathbb{R}^3$.
We study these links from several different
aspects including polynomial invariants using the relation between link polynomials of a transitive link and its factor links.
\end{abstract}


\maketitle

\section{Introduction}

Symmetry is one of the oldest and richest subjects not only in mathematics but also in many different disciplines
including engineering, designs, network models. Even though each discipline has different prospectives of symmetry,
a common interest is to realize the highest symmetry possible. In geometry and topology, symmetry plays a key role in
modern research. In the present article, we will study the links of the highest symmetry :``\emph{transitive links}".

A \emph{link} $L$ is an embedding of $n$ copies of $\mathbb{S}^1$ into $\mathbb{S}^3$. Since
we may consider $\mathbb{S}^3$ as $\mathbb{R}^3\cup \{\infty\}$, we will assume all links are in
$\mathbb{S}^3$ or $\mathbb{R}^3$ depending on our convenience.
If a link has only one copy of $\mathbb{S}^1$, the link is called a \emph{knot}. Two links are \emph{equivalent}
if there is an isotopy between them. In the case of prime knots, this equivalence is the same as the existence of an orientation preserving homeomorphism on $\mathbb{S}^3$,
which sends a knot to the other knot. Although the equivalent class of a link $L$ is called a \emph{link type}, throughout the article, a link
really means the equivalent class of link $L$. Additional terms in the knot theory can be found in~\cite{BZ:knot}.

One classical invariant in knot theory is the periodicity.
A link $L$ in $\mathbb{S}^3$ is $p$-\emph{periodic} if there exists an orientation preserving
periodic homeomorphism $h$ of order $p$ such that $fix(h)$ is homeomorphic to $S^1$,
$h(L)=L$ and $fix(h)\cap L=\emptyset$ where $fix(h)$ is the set of
fixed points of $h$. By the positive solution of Smith conjecture, $fix(h)$ is unknotted.
Thus, if we consider $\mathbb{S}^3$ as
$\mathbb{R}^3\cup \{\infty\}$, we can assume that $h$ is a rotation by
$2\pi/p$ angle around the $z-$axis. If $L$ is a periodic link, we
denote its factor link $(\mathbb{S}^3,L)/h$ by $\overline{L}$.
Murasugi~\cite{mu:alexander} found a strong relation between the
Alexander polynomials of a periodic link and its factor link.
Murasugi also found a similar relation for the Jones polynomials of
$L$ and $\overline L$~\cite{mu:jones}. There are various result to decide
periodicity of links \cite{JK, KL:sl3, Przytycki:criterion, traczyk:period3,
yokota:skein, yokota:kauffman}. These are all necessary conditions for periodic links using polynomial invariants of links.
There is no complete classification for periodic links yet.

For the periodicity of links,
the homeomorphism are all rotations. But, some of non periodic
links are invariant under some homeomorphisms which are not necessary rotations.
This motivates us to enlarge our interest for links of non rotational symmetries.

A \emph{symmetry group} of a link $L$ is the mapping class group of the pair $(\mathbb{S}^3, L)$.
More succinctly, a knot symmetry is a homeomorphism of the pair of spaces $(\mathbb{S}^3, K)$.
Hoste et al. \cite{HTJ} consider four types of symmetry based on whether the symmetry preserves
or reverses orienting of $\mathbb{S}^3$ and $K$, 1) preserves $\mathbb{S}^3$, preserves $K$
(identity operation), 2) preserves $\mathbb{S}^3$, reverses $K$, 3) reverses $\mathbb{S}^3$,
preserves $K$ and 4) reverses $\mathbb{S}^3$, reverses $K$. This then gives the five possible
classes of symmetry summarized in Table~\ref{table1}.

The symmetry groups of links have been studied very well in knot theory~\cite{BS, BZ, HW, KS}. Kodama and Sakuma used a
method in Bonahon and Siebenmann~\cite{BS} to compute these groups for all but three of the knots of $10$ and
fewer crossings~\cite{KS}. Henry and Weeks used the program SnapPea to compute the
symmetry groups for hyperbolic knots and links of $9$ and fewer crossings~\cite{HW}. These efforts followed
earlier tabulations of symmetry groups by Boileau and Zimmermann~\cite{BZ}, who found symmetry groups
for non-elliptic Montesinos links with $11$ or fewer crossings.
In the case of hyperbolic knots, the symmetry group must be finite and either cyclic or dihedral~\cite{HTJ, Riley, KS}.
The classification is slightly more complicated for nonhyperbolic knots. Furthermore,
all knots with $\le 8$ crossings are either amphichiral or invertible~\cite{HTJ}.
Any symmetry of a prime alternating link must be visible up to flypes
in any alternating diagram of the link~\cite{MT:alternating, HTJ}.
Hoste et al. found the numbers of $k$-crossing knots belonging to
cyclic symmetry groups and dihedral symmetry groups~\cite{HTJ}.

On the other hand, there is a very closely related but slight different
approach for the symmetry group of a link $L$, \emph{``intrinsic" symmetry group}.
Following the idea of Fox~\cite{Fox}, Whitten defined the group of symmetries of oriented,
labeled link $L$~\cite{Whitten}. J. Cantarella and et. al. find the intrinsic symmetry group
of links with $8$ and fewer crossings~\cite{BC:8}.

\begin{table}
\begin{tabular}{|c|c|c|}\hline
class  & symmetries & knot symmetries\\
  \hline
 $c$ & 1 & chiral, noninvertible \\ \hline
 $+$ & 1, 3 & $+$ amphichiral, noninvertible \\ \hline
 $-$ & 1, 4 & $-$ amphichiral, noninvertible \\ \hline
 $i$ & 1, 2 & chiral, invertible \\ \hline
 $a$ & 1, 2, 3, 4 &  $+$ and $-$ amphichiral, invertible
\\\hline
\end{tabular}
\vskip .5cm
\caption{Five types of the symmetry of knots} \label{table1}
\end{table}

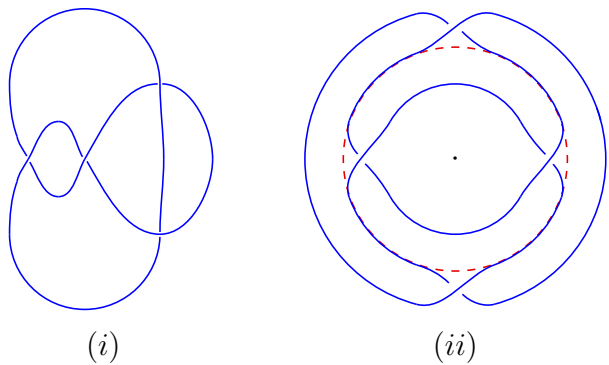
\begin{figure}
$$
\begin{pspicture}[shift=.2](-.2,-2.7)(2.9,2.2)
\psarc(1,1){1}{0}{180}
\pccurve[angleA=-90,angleB=110](0,1)(.1,.3)
\pccurve[angleA=-80,angleB=120](.1,.3)(.23,.04)
\pccurve[angleA=-60,angleB=180](.27,-.04)(.65,-.5)
\pccurve[angleA=0,angleB=-120](.65,-.5)(1,0)
\pccurve[angleA=60,angleB=180](1,0)(1.97,1)
\pccurve[angleA=0,angleB=90](2.03,1)(2.7,0)
\pccurve[angleA=-90,angleB=0](2.7,0)(2,-1)
\pccurve[angleA=180,angleB=-60](2,-1)(1.03,-.03)
\pccurve[angleA=120,angleB=0](.97,.03)(.65,.5)
\pccurve[angleA=180,angleB=60](.65,.5)(.25,0)
\pccurve[angleA=-120,angleB=80](.25,0)(.1,-.3)
\pccurve[angleA=-110,angleB=90](.1,-.3)(0,-1)
\psarc(1,-1){1}{180}{358}
\pccurve[angleA=90,angleB=-90](2,-.97)(2.05,0)
\pccurve[angleA=90,angleB=-90](2.05,0)(2,1)
\rput(1.25,-2.5){$(i)$}
\end{pspicture} \quad\quad
\begin{pspicture}[shift=.2](-2.2,-2.7)(2.2,2.2)
\psarc(0,0){2}{15}{75}
\pscircle[linecolor=darkred, linestyle=dashed](0,0){1.5}
\pccurve[angleA=165,angleB=20](2;75)(1.5;110)
\psarc(0,0){1.5}{110}{160}
\pccurve[angleA=-110,angleB=135](1.5;160)(1.3;177)
\pccurve[angleA=-45,angleB=120](1.2;183)(1;210)
\psarc(0,0){1}{210}{330}
\pccurve[angleA=60,angleB=-70](1;330)(1.5;20)
\psarc(0,0){1.5}{20}{70}
\pccurve[angleA=160,angleB=-45](1.5;70)(1.7;87)
\pccurve[angleA=135,angleB=15](1.8;93)(2;105)
\psarc(0,0){2}{105}{255}
\pccurve[angleA=-15,angleB=-160](2;255)(1.5;290)
\psarc(0,0){1.5}{290}{-20}
\pccurve[angleA=70,angleB=-45](1.5;-20)(1.3;-3)
\pccurve[angleA=135,angleB=-60](1.2;3)(1;30)
\psarc(0,0){1}{30}{150}
\pccurve[angleA=-120,angleB=110](1;150)(1.5;200)
\psarc(0,0){1.5}{200}{250}
\pccurve[angleA=-15,angleB=135](1.5;250)(1.65;267)
\pccurve[angleA=-45,angleB=195](1.8;273)(2;285)
\psarc(0,0){2}{285}{20}
\rput(0,0){$\cdot$}
\rput(0,-2.5){$(ii)$}
\end{pspicture}
$$
\caption{$(i)$ A standard diagram and $(ii)$ a transitive diagram of the figure eight knot.}
\label{Figure8}
\end{figure}

The origin of the present article is the symmetry group of the figure eight knot.
It is already known that the symmetry group of the figure eight knot is the dihedral group
of order $8$ which not only contains two obvious horizontal and vertical reflections,
but also contains an element of order $4$, the composition of the reflection along the the dashed red circle
and the rotation by $\dfrac{\pi}{2}$ along the point as illustrated in Figure~\ref{Figure8}.
This element of order $4$ transitively acts on the set of crossings. This phenomenon naturally raises a new direction of
the study of symmetry of links. A link diagram $D(L)$ of a link is \emph{transitive} if the symmetry group of the link $L$
acts transitively on the set of crossings in the diagram. A link $L$ is \emph{transitive} if it admits a transitive diagram.
The benefit of having transitive diagram is numerous; first we can import some known results in graph theory to find
a complete classification of transitive links, second we can extend this classification for periodic link and we can obtained
several necessary conditions for being a transitive link using link polynomials.

\begin{reptheorem}{transitiveS3}
A link $L$ is transitive if and only if it is either $(2n+1)_1$, $(2n)_1^2$ in Figure~\ref{Cn},
$Ch_n$ in Figure~\ref{RP}, $\overline{(\sigma_1 \sigma_2^{-1})^n}$ in Figure~\ref{RAP}
or one of the eight links corresponding to the eight Archimedian solids
depicted in Figure~\ref{TruncatedTetrahedron} to Figure~\ref{TruncatedDodecahedron}.
\end{reptheorem}

The outline of this paper is as follows. We first provide some preliminary definitions and related results in graph theory
 in Section~\ref{prelim}. In Section\ref{transitivelink},
we investigate the transitive links and find a complete classification of them including the proof of Theorem~\ref{transitiveS3}.
In Section~\ref{disscusion}, we conclude with further research problems.

\section{Preliminaries} \label{prelim}

Let us give a list of definitions we will be using throughout the rest of article.

Let $L$ be a link. Let $\Sym(L)$ be the set of all homeomorphisms on $\mathbb{S}^3$
which preserves $L$. Let $\Sym^+(L)$ be the set of all orientation preserving homeomorphisms on $\mathbb{S}^3$ which preserves $L$.
\begin{defi}
A link $L$ is \emph{transitive} if it admits a link diagram on a plane in $\mathbb{R}^3$
such that its $\Sym(L)$ acts transitively on the set of all crossings of $L$.
A link $L$ is \emph{positive transitive} if it is transitive with respect to $\Sym^+(L)$.

A link $L$ is \emph{tangle transitive} if there exists a nontrivial, which means the numbers of tangles must be at least $2$,
tangle decomposition $\mathcal{T}(L)=\{T_1, T_2, \ldots, T_n\}$ of $L$ on a
plane in $\mathbb{R}^3$ and $\Sym(L)$ acts transitively on $\mathcal{T}(L)$.
A link $L$ \emph{positive tangle transitive} if it is tangle transitive with respect to $\Sym^+(L)$.

A link $L$ \emph{block transitive} if there exists a nontrivial (the numbers of blocks must be at least $2$)
blocks decomposition $\mathcal{B}(L)=\{B_1, B_2, \ldots, B_n\}$ of $L$ on a plane in
$\mathbb{R}^3$and $\Sym(L)$ acts transitively on $\mathcal{B}(K)$.
A link $L$ \emph{positive tangle transitive} if it is block transitive with respect to $\Sym^+(L)$.
For example, the figure eight knot is transitive but not positive transitive.
\end{defi}

A \emph{graph} $\Gamma$ is an ordered pair $\Gamma = (V(\Gamma), E(\Gamma))$ comprising a set
$V(\Gamma)$ of vertices together with a set $E(\Gamma)$ of edges. Two graphs $\Gamma_1 = (V(\Gamma_1), E(\Gamma_1))$
and $\Gamma_2 = (V(\Gamma_2), E(\Gamma_2))$ are \emph{equivalent} if there exists a bijective function $\phi :  V(\Gamma_1) \longrightarrow V(\Gamma_2)$
such that $e =\{ u, v \} \in E(\Gamma_1)$ if and only if $\{ \phi(u), \phi(v) \} \in E(\Gamma_2)$ and $\phi$ is called a \emph{graph isomorphism}.
If $\Gamma_1 =\Gamma_2$, the graph isomorphism $\phi$ is often called a graph \emph{automorphism}. The set of all graph automorphism on $\Gamma$ is
\emph{graph automorphism group} denoted by $\Aut(\Gamma)$.

Not surprisingly, the transitivity is not new in graph theory.
\begin{defi}
A graph $\Gamma$ is \emph{vertex transitive} if $\Aut(\Gamma)$ acts transitively
on the set of vertices $V(\Gamma)$. A graph $\Gamma$ is \emph{edge transitive} if $\Aut(\Gamma)$ acts transitively
on $E(\Gamma)$. A graph $\Gamma$ is \emph{arc transitive} or \emph{symmetric} if for
$e_1 =(u_1, v_1), e_2=(u_2, v_2) \in E(\Gamma)$ there exists $ \phi \in\Aut(\Gamma)$ such that
$\phi(u_1)=u_2$ and $\phi(v_1)=v_2$.
\end{defi}

Every symmetric graph without isolated vertices is vertex and edge transitive,
and every vertex-transitive graph is regular.
However, not all vertex-transitive graphs are symmetric, for example, the edges of the truncated tetrahedron,
and not all regular graphs are vertex-transitive, for example, the Frucht graph and Tietze's graph.
There have been a serious study on vertex transitive, edge and arc transitive graphs[Add Reference].

A typical example of vertex transitive graphs is a Cayley graph.
\begin{defi}
Let $G$ be a group and $S$ be a subset
of $G$ which is called a \emph{generating set}.
\emph{The Cayley graph $\Gamma=(G,S)$} is colored directed graph whose
vertex set $V(\Gamma)=G$, for each generator $s \in S$ is assigned a color $c_s$, and the edge set
$E(\Gamma)=\{ (g,gs) | g \in G, s \in S \}$ where the edge $(g,gs)$ is colored by $c_s$.
\end{defi}

Not all vertex-transitive graphs are Cayley graphs, for example, the peterson graph is vertex transitive but is not a Cayley graph.
For our purpose to relate with transitive links, these vertex transitive graphs have to be planar. It is already known that
vertex transitive simple graphs of valency $>5$ is not planar.

There have been numerous results about the transitive graphs and Cayley graphs.
The following two theorems will be used for our classification theorems in Section~\ref{transitivelink}.

\begin{theorem} (\cite{Maschke:cayley}) \label{Maschkethm}
The only groups that can give planar Cayley graphs are exactly $\mathbb{Z}_n, \mathbb{Z}_2 \times \mathbb{Z}_n, D_{2n}, S_4, A_4$ and $A_5$.
\end{theorem}

\begin{theorem} (\cite{FI:transitive}) \label{FIthm}
A connected simple graph $\Gamma$ is planar vertex transitive graphs if and only if it is either
a point, $K_2$, $C_n$, regular prisms, regular anti-prisms, the Platonic solids or the Archimedean solids.
\end{theorem}

\section{Transitive links} \label{transitivelink}

A key observation to find a complete classification of the transitive links is that
if we ignore the crossing of a link diagram to make it double point, it become a planar graph of valency $4$.
Furthermore, if the link diagram is transitive, the resulting graph is vertex transitive.
The converse is not true in general. For example, the octahedron is a planar vertex transitive graph but
it is also edge transitive. Thus, even if we fix one of crossing types at a vertex, by the action of
the graph automorphism group of the octahedron we may not recover a link diagram.
This observation leads us to the following lemma.

\begin{lem} \label{val4}
A planar vertex transitive graph $\Gamma$ of valency $4$ is obtained from a transitive link diagram by projecting crossings to double points
if and only if $\Gamma$ is not edge transitive.
\begin{proof}

\end{proof}
\end{lem}

By combining Lemma~\ref{val4} and Theorem~\ref{Maschkethm}, we find the following theorem.

\begin{theorem}
Planar Cayley diagrams of valency $4$ are the truncated tetrahedron, the cubeoctahedron, the truncated octahedron,
the truncated cube, the rhombicuboctahedron, the truncated icosahedron, the truncated dodecahedron and the rhombicosidodecahedron.
\end{theorem}

\begin{figure}
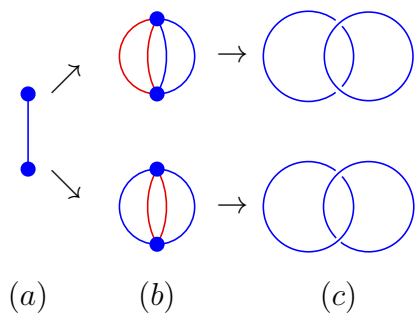

$$

$$
\caption{$(a)$ The complete graph $K_2$, $(b)$ its possible multi-graphs and $(c)$ corresponding transitive links.}
\label{K2}
\end{figure}

\begin{figure}
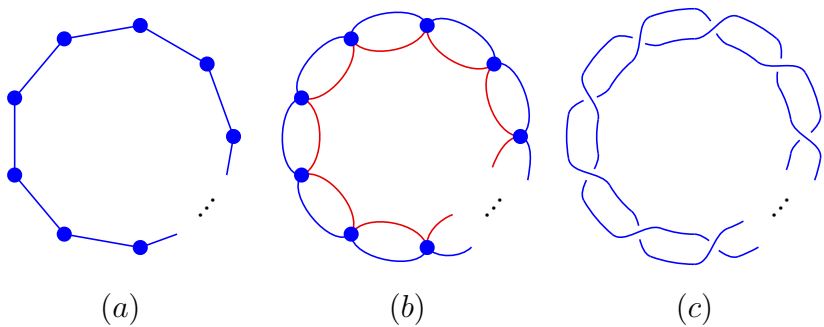

$$
%
$$
\caption{$(a)$ The cycle graph $C_n$ of $n$ vertices, $(b)$ its possible multi-graph of valency $4$ and $(c)$ corresponding transitive link $\overline{(\sigma_1)^n}$.}
\label{Cn}
\end{figure}

\begin{figure}
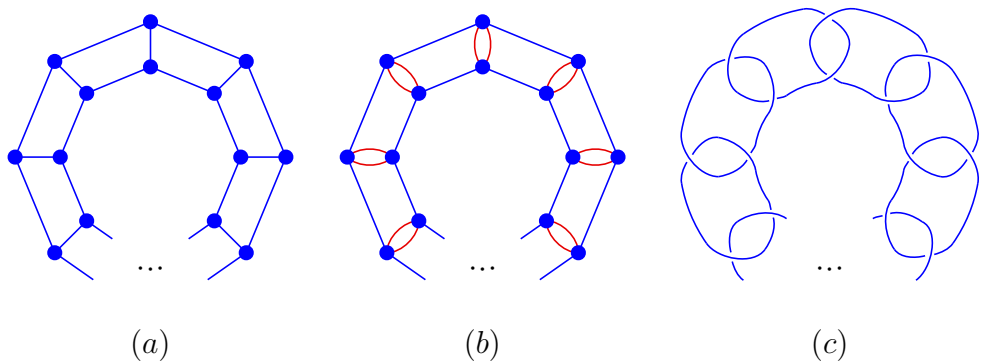

$$
%
  $$
\caption{$(a)$ The regular prism, $(b)$ its possible multi-graph of valency $4$ and $(c)$ corresponding transitive links.}
\label{RP}
\end{figure}

\begin{figure}
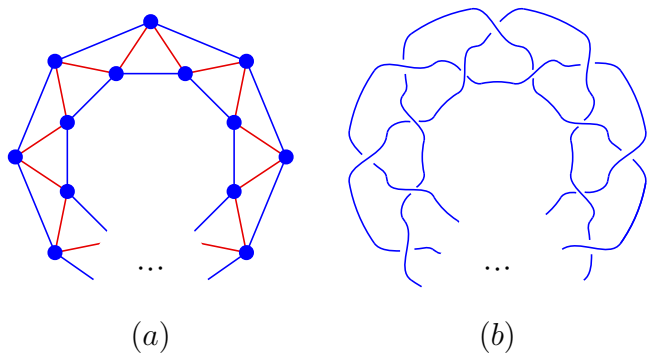

$$
%
  $$
\caption{$(a)$ The regular antiprism and $(b)$ its corresponding transitive link $\overline{(\sigma_1 \sigma_2^{-1})^n}$.}
\label{RAP}
\end{figure}

\begin{figure}
$$
%
 $$
\caption{$(a)$ The Truncated Tetrahedron, $(b)$ its possible multi-graph of valency $4$ and $(c)$ corresponding transitive link.}
\label{TruncatedTetrahedron}
\end{figure}

\begin{figure}
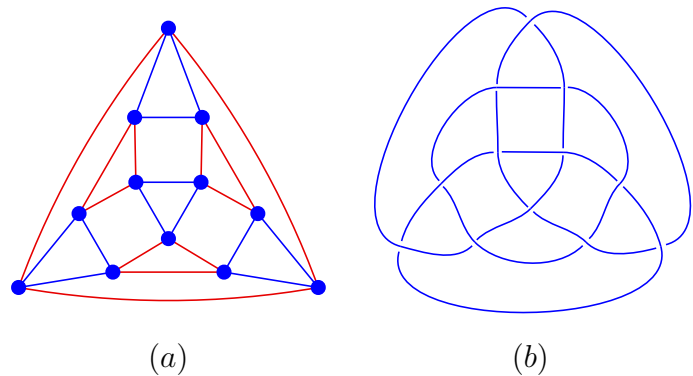

$$
%
 $$
\caption{$(a)$ The Cube Octahedron and $(b)$ its corresponding transitive link.}
\label{CubeOctahedron}
\end{figure}

\begin{figure}
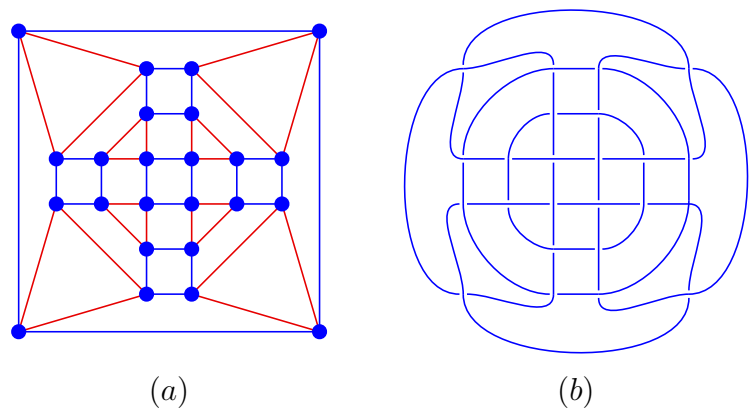

$$
%
 $$
\caption{$(a)$ The Rhombicuboctahedron and $(b)$ its corresponding transitive link.}
\label{Rhombicuboctahedron}
\end{figure}

\begin{figure}
$$
 %
 $$
\caption{$(a)$ The Truncated Octahedron, $(b)$ its possible multi-graph of valency $4$ and $(c)$ corresponding transitive link.}
\label{TruncatedOctahedron}
\end{figure}

\begin{figure}
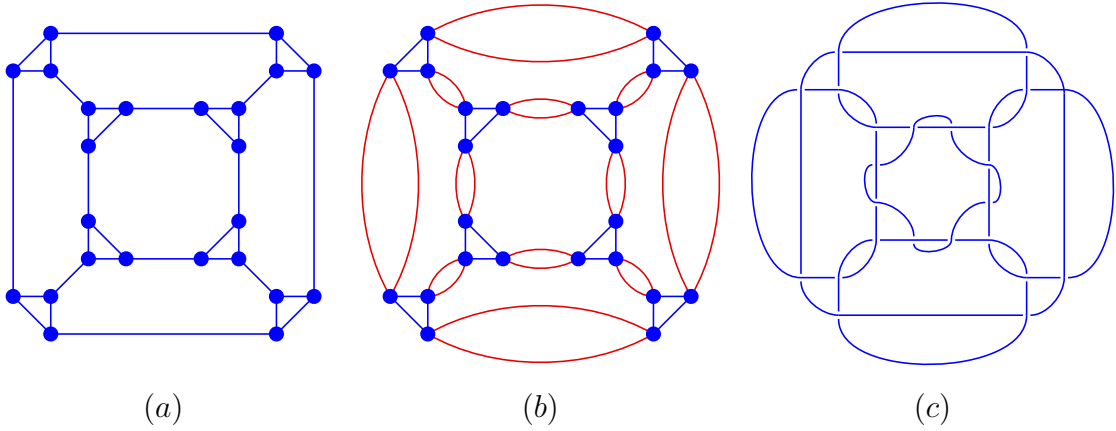

$$
 %
 $$
\caption{$(a)$ The Truncated Cube, $(b)$ its possible multi-graph of valency $4$ and $(c)$ corresponding transitive link.}
\label{TruncatedCube}
\end{figure}

\begin{figure}
$$
%
 $$
\caption{$(a)$ The Truncated icosahedron, $(b)$ its possible multi-graph of valency $4$ and $(c)$ corresponding transitive link.}
\label{TruncatedIcosahedron}
\end{figure}

\begin{figure}
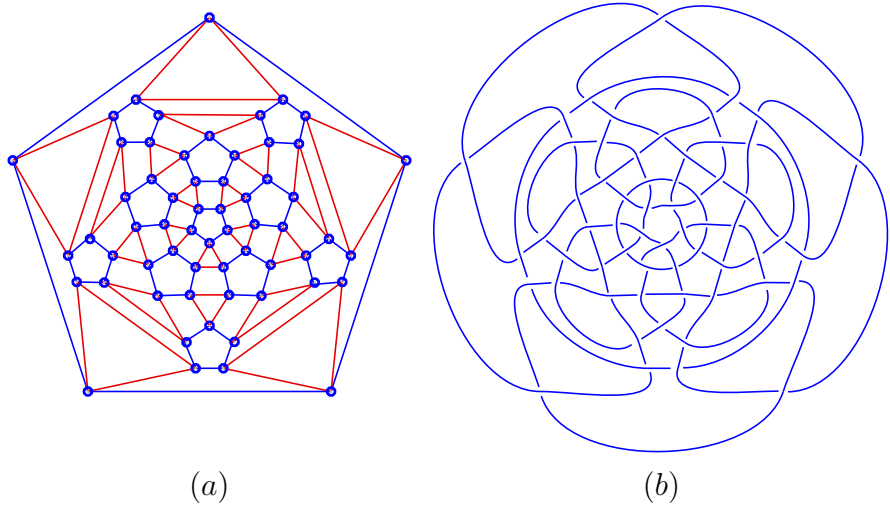

$$
%
 $$
\caption{$(a)$ The Rhombicosidodecahedron and $(b)$ its corresponding transitive link.}
\label{Rhombicosidodecahedron}
\end{figure}

\begin{figure}
$$
%
  $$
\caption{$(a)$ The Truncated Dodecahedron, $(b)$ its possible multi-graph of valency $4$ and $(c)$ corresponding transitive link.}
\label{TruncatedDodecahedron}
\end{figure}

By summarizing these, I can obtain the following theorem.

\begin{theorem} \label{transitiveS3}
A link $L$ is transitive if and only if it is either $(2n+1)_1$, $(2n)_1^2$ in Figure~\ref{Cn},
$Ch_n$ in Figure~\ref{RP}, $\overline{(\sigma_1 \sigma_2^{-1})^n}$ in Figure~\ref{RAP}
or one of the eight links corresponding to the eight Archimedian solids
depicted in Figure~\ref{TruncatedTetrahedron} to Figure~\ref{TruncatedDodecahedron}.
\begin{proof}
From the single vertex by adding two loops, we get the bouquet of $2$ circle.
No matter how we put a crossing, we get the unknot.
From the single edge by adding three edges, we get the dipole graph $D_4$.
Once we fix a crossing at a vertex,
By the rotation we get the Hope link and by the reflection, we get the trivial links of two components.
From the simple circuit by making all double edges, we get the $\overline{(\sigma_1)^n}$ on two strings.
For the regular prism $\Gamma$ which has valency $3$, the orbit of
edges by $Aut(\Gamma)$ are two, the edges joining
top and bottom regular polygon can be replaced by double edges,
we get the following link $(4n)_1^{n}$ as in Figure~\ref{RP}.
From the regular anti-prism which has valency $4$, by the action
of $\mathbb{Z}_2 \times \mathbb{Z}_n$ we get $\overline{(\sigma_1 \sigma_2)^{\pm n}}$,
by the action of $D_{2n}$ we get $\overline{(\sigma_1 \sigma_2^{-1})^{\pm n}}$ on three strings.
For the Platonic solids $\Gamma$ which has valency $3$ is also edge transitive,
so if I replace an edge by a double edge to make valency $4$, by the action of $Aut(\Gamma)$
it becomes a graph of valency $6$. So we rule out all the Platonic solids except the
cube which may be considered as a rectangular prism and octahedron which may be considered as a triangular anti-prism.
For the Archimedian solids, if we rule out Archimedian solids of valency $5$,
all remaining $8$ graphs are indeed Cayley graph of valency $4$ which are the truncated tetrahedron,
the cubeoctahedron, the truncated octahedron,
the truncated cube, the rhombicuboctahedron, the truncated icosahedron,
the truncated dodecahedron and the rhombicosidodecahedron.
If we fix a crossing at a vertex, by the action of the corresponding Cayley group,
we get $8$ links depicted in Figure~\ref{TruncatedTetrahedron} to Figure~\ref{TruncatedDodecahedron}.
\end{proof}
\end{theorem}

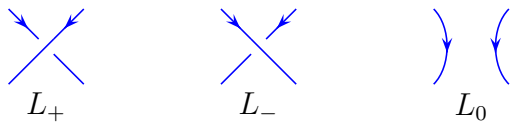
\begin{figure}$$
\begin{pspicture}[shift=-.8](0,-1)(0,.7) \end{pspicture} \begin{pspicture}[shift=-.8](0,.8)(0,.8)
\begin{pspicture}[shift=-.8](-1,-1)(1,.7) \rput(.5,.5){\rnode{a1}{$$}} \rput(-.5,.5){\rnode{a2}{$$}}
\rput(-.5,-.5){\rnode{a3}{$$}} \rput(.5,-.5){\rnode{a4}{$$}} \rput(.1,.1){\rnode{b1}{$$}} \rput(-.1,.1){\rnode{b2}{$$}}
\rput(-.1,-.1){\rnode{b3}{$$}} \rput(.1,-.1){\rnode{b4}{$$}} \ncline{a1}{b1}\middlearrow  \ncline{a2}{b2}\middlearrow
\ncline{b4}{a4} \ncline{b1}{a3} \rput(0,-.8){\rnode{c4}{$L_+$}} \end{pspicture}\end{pspicture} \quad \quad
\begin{pspicture}[shift=-.8](0,.8)(0,.8) \begin{pspicture}[shift=-.8](-1,-1)(1,.7) \rput(.5,.5){\rnode{a1}{$$}} \rput(-.5,.5){\rnode{a2}{$$}}
\rput(-.5,-.5){\rnode{a3}{$$}} \rput(.5,-.5){\rnode{a4}{$$}} \rput(.1,.1){\rnode{b1}{$$}} \rput(-.1,.1){\rnode{b2}{$$}} \rput(-.1,-.1){\rnode{b3}{$$}}
\rput(.1,-.1){\rnode{b4}{$$}} \ncline{a2}{b2}\middlearrow \ncline{a1}{b1}\middlearrow \ncline{b2}{a4} \ncline{b3}{a3}
\rput(0,-.8){\rnode{c4}{$L_-$}} \end{pspicture}\end{pspicture}
\quad  \quad \begin{pspicture}[shift=-.8](0,.8)(0,.8) \begin{pspicture}[shift=-.8](-1,-1)(1,.7)
\rput(.5,.5){\rnode{a1}{$$}} \rput(-.5,.5){\rnode{a2}{$$}}
\rput(-.5,-.5){\rnode{a3}{$$}} \rput(.5,-.5){\rnode{a4}{$$}}
\rput(.1,.1){\rnode{b1}{$$}} \rput(-.1,.1){\rnode{b2}{$$}}
\rput(-.1,-.1){\rnode{b3}{$$}} \rput(.1,-.1){\rnode{b4}{$$}}
\nccurve[angleA=225,angleB=135]{a1}{a4}\middlearrow
\nccurve[angleA=315,angleB=45]{a2}{a3}\middlearrow
\rput(0,-.8){\rnode{c4}{$L_0$}}
\end{pspicture}\end{pspicture}$$
\caption{The skein triple $L_+, L_-$ and $L_0$.} \label{local}
\end{figure}

Now we define the HOMFLY polynomial specialized to a one variable polynomial.
For a nonnegative integer $n$, the HOMFLY polynomial $P_n(q)$ specialized to a one variable polynomial can be calculated
uniquely by the following skein relations:

\begin{gather}\nonumber P_{n}(\emptyset) =1, \\\nonumber P_{n}(
\begin{pspicture}[shift=-.07](-.17,-.17)(.17,.17) \pscircle(0,0){.15}
\end{pspicture} \cup D)= (\frac{q^{\frac n2}- q^{-\frac n2}}{q^{\frac 12}- q^{-\frac 12}}) P_{n}(D),
 \\\nonumber q^{\frac n2}P_{n}(L_+) - q^{-\frac n2}P_{n}(L_-) = (q^{\frac 12}- q^{-\frac 12}) P_n(L_0), \end{gather}
where $\emptyset$ is the empty diagram,
$\begin{pspicture}[shift=-.07](-.17,-.17)(.17,.17) \pscircle(0,0){.15} \end{pspicture}$
is the trivial knot and $L_+, L_-$ and $L_0$ are skein triple, three diagrams which are identical except at
one crossing as in Figure~\ref{local}.

\begin{theorem} (\cite{JK})
Let $p$ be a positive integer and $L$ be a $p$-periodic link in $\mathbb{S}^3$
with its factor link $\overline{L}$. Then,
$$P_n(L) \equiv P_n(\overline{L})^p  \hskip .7cm modulo \hskip .2cm
\mathcal{I}_n,$$ where $\mathcal{I}_n$ is the ideal of
$\mathbb{Z}[q^{\pm \frac 12}]$ generated by $p$ and
$\left[\begin{matrix}n\\
i\end{matrix}\right]^p-\left[\begin{matrix}n\\
i\end{matrix}\right]$ for $i=1, 2, \ldots, \lfloor \frac{n}{2}
\rfloor$. \label{modifiedconj}
\end{theorem}

A precise and algebraic overview of the quantum $\mathfrak{sl}(n)$ representation theory can be found in ~\cite{Morrison}.
If links are decorated by the fundamental representations
$V_{\lambda_i}$ of the quantum $\mathfrak{sl}(n)$, denoted by $i$,
Murakami, Ohtsuki and Yamada~\cite{MOY:Homfly} found a quantum
invariant $[D]_n$ for framed links by resolving each crossing in a link
diagram $D$ of $L$ as shown in Equation~(\ref{pstvresol1}) and
Equation~(\ref{pstvresol2}) in Figure~\ref{expansion}.
For a coloring $\mu$ (a representation of the quantum $\mathfrak{sl}(n)$) of a diagram $D$ of a
link $L$, we first consider a \emph{colored writhe} $\omega_i(D)$ as
the sum of writhes of components colored by $i$. Then we set
$$K_n(L,\mu) = \prod_i q^{-\omega_i(D)\frac{i(n-i+1)}{2}} [D]_n,$$
where the product runs over all colors $i$.

\begin{figure}
\begin{eqnarray}\left[
\begin{pspicture}[shift=-.6](-.5,-.7)(.5,.7) \rput[bl](.3,.6){$j$}
\rput[br](-.3,.6){$i$} \rput[tr](-.3,-.6){$j$}
\rput[tl](.3,-.6){$i$} \psline(-.25,-.5)(.25,.5)
\psline[arrowscale=1.5]{->}(.15,.3)(.25,.5) \psline(.25,-.5)(.05,-.1)
\psline[arrowscale=1.5]{->}(-.15,.3)(-.25,.5) \psline(-.05,.1)(-.25,.5)
\end{pspicture}\right]_n
=\sum_{k=0}^{i} (-1)^{k+(j+1)i}q^{\frac{(i-k)}{2}} \left[
\begin{pspicture}[shift=-1.2](-1.75,-1.2)(1.75,1.2)
\rput[bl](1.2,.85){$j$} \rput[br](-1.2,.85){$i$}
\rput[tr](-1.2,-.85){$j$} \rput[tl](1.2,-.85){$i$}
\rput[b](0,.5){$j+k-i$} \rput[r](-.9,0){$j+k$}
\rput[t](0,-.5){$k$} \rput[l](.9,0){$i-k$}
\psline(.8,.4)(1.2,.8)\psline[arrowscale=1.5]{->}(.9,.5)(1.1,.7)
\psline(-.8,.4)(-1.2,.8)\psline[arrowscale=1.5]{->}(-.9,.5)(-1.1,.7)
\psline(-1.2,-.8)(-.8,-.4)\psline[arrowscale=1.5]{->}(-1.1,-.7)(-.9,-.5)
\psline(1.2,-.8)(.8,-.4)\psline[arrowscale=1.5]{->}(1.1,-.7)(.9,-.5)
\psline(-.8,.4)(.8,.4)\psline[arrowscale=1.5]{->}(-.1,.4)(.1,.4)
\psline(.8,-.4)(.8,.4)\psline[arrowscale=1.5]{->}(.8,-.1)(.8,.1)
\psline(-.8,-.4)(-.8,.4)\psline[arrowscale=1.5]{->}(-.8,-.1)(-.8,.1)
\psline(.8,-.4)(-.8,-.4)\psline[arrowscale=1.5]{->}(.1,-.4)(-.1,-.4)
\end{pspicture}\right]_n  \label{pstvresol1}
\end{eqnarray}
\begin{eqnarray}\left[
\begin{pspicture}[shift=-.6](-.5,-.7)(.5,.7) \rput[bl](.3,.6){$j$}
\rput[br](-.3,.6){$i$} \rput[tr](-.3,-.6){$j$}
\rput[tl](.3,-.6){$i$} \psline(-.25,-.5)(.25,.5)
\psline[arrowscale=1.5]{->}(.15,.3)(.25,.5) \psline(.25,-.5)(.05,-.1)
\psline[arrowscale=1.5]{->}(-.15,.3)(-.25,.5) \psline(-.05,.1)(-.25,.5)
\end{pspicture}\right]_n
=\sum_{k=0}^{j} (-1)^{k+(i+1)j}q^{\frac{(j-k)}{2}}\left[
\begin{pspicture}[shift=-1.2](-1.75,-1.2)(1.75,1.2)
\rput[bl](1.2,.85){$j$} \rput[br](-1.2,.85){$i$}
\rput[tr](-1.2,-.85){$j$} \rput[tl](1.2,-.85){$i$}
\rput[b](0,.5){$i+k-j$} \rput[r](-.9,0){$j-k$}
\rput[t](0,-.5){$k$} \rput[l](.9,0){$i+k$}
\psline(.8,.4)(1.2,.8)\psline[arrowscale=1.5]{->}(.9,.5)(1.1,.7)
\psline(-.8,.4)(-1.2,.8)\psline[arrowscale=1.5]{->}(-.9,.5)(-1.1,.7)
\psline(-1.2,-.8)(-.8,-.4)\psline[arrowscale=1.5]{->}(-1.1,-.7)(-.9,-.5)
\psline(1.2,-.8)(.8,-.4)\psline[arrowscale=1.5]{->}(1.1,-.7)(.9,-.5)
\psline(-.8,.4)(.8,.4)\psline[arrowscale=1.5]{<-}(-.1,.4)(.1,.4)
\psline(.8,-.4)(.8,.4)\psline[arrowscale=1.5]{->}(.8,-.1)(.8,.1)
\psline(-.8,-.4)(-.8,.4)\psline[arrowscale=1.5]{->}(-.8,-.1)(-.8,.1)
\psline(.8,-.4)(-.8,-.4)\psline[arrowscale=1.5]{<-}(.1,-.4)(-.1,-.4)
\end{pspicture}\right]_n  \label{pstvresol2}
\end{eqnarray}
\caption{Skein expansions of a crossing} \label{expansion}
\end{figure}

By decorating each component by $\mu$,
we can define the colored $\mathfrak{sl}(n)$
HOMFLY polynomial, denoted by $G_n(L,\mu)$ as follows.
 For a given colored link $L$ of $l$ components
say, $L_1, L_2, \ldots, L_l$, where each component $L_i$ is colored
by an irreducible representation $V_{a(i)_1\lambda_1 +
a(i)_2\lambda_2 + \ldots + a(i)_{n-1}\lambda_{n-1}}$ of
$\mathfrak{sl}(n)$ and $\lambda_1$, $\lambda_2$, $\ldots$, $\lambda_{n-1}$
are the fundamental weights of $\mathfrak{sl}(n)$. The coloring is
denoted by $\mu=(a(1)_1\lambda_1 + a(1)_2\lambda_2 + \ldots +
a(1)_{n-1}\lambda_{n-1}, a(2)_1\lambda_1 + a(2)_2\lambda_2 + \ldots
+ a(2)_{n-1}\lambda_{n-1}, \ldots, a(l)_1\lambda_1 + a(l)_2\lambda_2
+ \ldots + a(l)_{n-1}\lambda_{n-1})$. First we replace each
component $L_i$ by $a(i)_1+a(i)_2+\ldots +a(i)_{n-1}$ copies of
parallel lines and each $a(i)_j$ line is colored by the weight
$\lambda_j$. Then we put a clasp of weight $(a(i)_1\lambda_1 +
a(i)_2\lambda_2 + \ldots + a(i)_{n-1}\lambda_{n-1})$ for $L_i$. If
we assume the clasps are far away from crossings, we expand each
crossing as in Figure~\ref{expansion}, then
clasps \cite{JK, yokota:skeinforn}. The value we obtained after removing all
faces by using the relations is \emph{the colored $\mathfrak{sl}(n)$
HOMFLY polynomial} $G_n(L,\mu)$ of $L$.

\begin{theorem} (\cite{JK})
Let $p$ be a positive integer and $L$ be a $p-$periodic link in
$\mathbb{S}^3$ with its factor link $\overline{L}$. Let $\mu$ be a
$p$-periodic coloring of $L$ and $\overline{\mu}$ be the induced
coloring of $\overline{L}$. Then for $n\ge 0$,

$$K_n(L,\mu) \equiv K_n(\overline{L},\overline{\mu})^p \hskip .7cm modulo \hskip .2cm
\mathcal{I}_n,$$ where $\mathcal{I}_n$ is the ideal of
$\mathbb{Z}[q^{\pm \frac 12}]$ generated by $p$ and
$\left[\begin{matrix}n\\
i\end{matrix}\right]^p-\left[\begin{matrix}n\\
i\end{matrix}\right]$ for $i=1, 2, \ldots, \lfloor \frac{n}{2}
\rfloor$. \label{mainthm}
\end{theorem}

\begin{theorem} (\cite{JK})
Let $p$ be a positive integer and $L$ be a $p-$periodic link in
$\mathbb{S}^3$ with its factor link $\overline{L}$. Let $\mu$ be a
$p$-periodic coloring of $L$ and $\overline{\mu}$ be the induced
coloring of $\overline{L}$. Then for $n\ge 0$,

$$G_n(L,\mu) \equiv G_n(\overline{L},\overline{\mu})^p \hskip .7cm modulo \hskip .2cm
\mathcal{I}_n,$$ where $\mathcal{I}_n$ is the ideal of
$\mathbb{Z}[q^{\pm \frac 12}]$ generated by $p$ and
$\left[\begin{matrix}n\\
i\end{matrix}\right]^p-\left[\begin{matrix}n\\
i\end{matrix}\right]$ for $i=1, 2, \ldots, \lfloor \frac{n}{2}
\rfloor$. \label{mainthm4}
\end{theorem}

Using the relation between HOMFLY polynomials $P_n(q)$ specialized to a one variable polynomial as stated in Theorem~\ref{modifiedconj} and
the (colored, resp.) $\mathfrak{sl}(n)$ HOMFLY polynomial $K_n(*,\mu)$ ($G_n(*,\mu)$, resp.) specialized to a one variable polynomial of a
periodic link $L$ and its factor link $\overline L$ as stated in Theorem~\ref{mainthm} and Theorem~\ref{mainthm4}, we find the following necessary condition of being a transitive link because the factor
link is the unknot.

\begin{theorem} \label{nilpo}
Let $\mu$ be an irreducible representation of the quantum Lie algebra $\mathfrak{sl}(n)$ and
$(O, \mu)$ be the unknot colored by $\mu$.
For a positive integer $m$, let $\mathcal{I}_n$ be the ideal of
$\mathbb{Z}[q^{\pm \frac 12}]$ generated by $m$ and
$\left[\begin{matrix}n\\
i\end{matrix}\right]^p-\left[\begin{matrix}n\\
i\end{matrix}\right]$ for $i=1, 2, \ldots, \lfloor \frac{n}{2}
\rfloor$.
If a link $L$ is $m$-transitive, then
\begin{enumerate}
\item[{\rm (1)}] $P_n(L)=[n]^m$ modulo $\mathcal{I}_n$.

\item[{\rm (2)}] $K_n(L,\mu)=(K_n(O,\mu))^m$ modulo $\mathcal{I}_n$.

\item[{\rm (3)}] $G_n(L,\mu)=(G_n(O,\mu))^m$ modulo $\mathcal{I}_n$.

\end{enumerate}
\end{theorem}

Let us remark that if we replace the convention $P_{n}(\emptyset) =1$ by $P_{n}(
\begin{pspicture}[shift=-.07](-.17,-.17)(.17,.17) \pscircle(0,0){.15}
\end{pspicture}) =1$, then the statement in Theorem~\ref{nilpo} {\rm (1)}
can be restate that $P_n(L)$ is nilpotent modulo $\mathcal{I}_n$.

\section{Conclusion} \label{disscusion}

One may found similar classification theorems of vertex transitive graphs on torus and projective plane by Carsten Thomassen~\cite{Thomassen}.
If one wants to extend our results for these vertex transitive graphs, one might
have to choose transitive links in right $3$-manifold (maybe torus$\times I$ or torus$\times \mathbb{S}^1$), but
we do not know yet. We believe the problem is the symmetry group for
transitive link (diagram on $\mathbb{S}^2$) or in $\mathbb{S}^3$ are just not seriously different from $Aut(\Gamma)$ where
$\Gamma$ is the crossingless graph of the transitive link while this phenomena no longer
works for the vertex transitive graphs on torus and projective plane.

\vskip 1cm

 \noindent{\bf Acknowledgements}
The \TeX\, macro package PSTricks~\cite{PSTricks} was essential for
typesetting the equations and figures. This work was supported by Kyonggi University Research Grant 2012.


\begin{thebibliography}{00}

\bibitem{BC:8} M. Berglund, J. Cantarella, M. P. Casey, E. Dannenberg,
W. George, A. Johnson, A. Kelley, A. LaPointe, M. Mastin, J. Parsley,
J. Rooney and R. Whitaker, \textit{Intrinsic symmetry groups of links with 8 and fewer crossings}, Symmetry 4 (2012), 143--207.

\bibitem{BS} F. Bonahon, L. Siebenmann, \textit{New Geometric Splittings of Classical Knots, and the
Classification and Symmetries of Arborescent Knots}, preprint. Available at {\tt http://www-bcf.usc.edu/\~{ }fbonahon/
Research/Preprints/Preprints.html}

\bibitem{BZ:knot} G. Burde and H. Zieschang, Knots, Berlin, de Gruyter, 1985.

\bibitem{BZ} M. Boileau, B. Zimmermann, \textit{Symmetries of nonelliptic Montesinos links}, Math. Ann. 277 (1987),
563--584.

\bibitem{FI:transitive} H. Fleischner and W. Imrich, \textit{Transitive planar graphs}, Math. Slovaca 29 (1979), 97--106.

\bibitem{Fox} R. Fox, \textit{Some Problems in Knot Theory}, In Topology of 3-Manifolds And Related Topics
(Proc. The Univ. of Georgia Institute, 1961) Prentice-Hall, Englewood Cliffs, NJ, USA, 1962, pp. 168--176.

\bibitem{HW} S. Henry, J. Weeks,  \textit{Symmetry groups of hyperbolic knots and links}, J. Knot Theory Ramif.
1 (1992), 185--201.

\bibitem{HTJ} J. Hoste, M. Thistlethwaite and J. Weeks, \textit{The First $1701936$  Knots}, Math. Intell. 20 (1998), 33--48.

\bibitem{JK} M. Jeong and D. Kim, \textit{The quantum $\mathfrak{sl}(n,\mathbb{C})$ representation theory and its applications},
H. Korean Math. Soc. 49 (2012), 993--1015.

\bibitem{Jones:braid} V. F. R. Jones, \textit{Hecke algebra representations of braid groups and link polynomials},
Ann. of Math. 126 (1987), 335--388.


\bibitem{KL:sl3} D. Kim and J. Lee, \textit{The quantum sl(3) invariants of cubic bipartite planar graphs},
J. Knot Theory Ramifications, 17(3) (2008), 361--375.


\bibitem{KS} K. Kodama and M. Sakuma, \textit{Symmetry Groups of Prime Knots Up to $10$ Crossings}, In Knot 90,
Proceedings of the International Conference on Knot Theory and Related Topics, Osaka, Japan, 1990 (Ed. A. Kawauchi.) Berlin: de Gruyter, pp. 323--340, 1992.

\bibitem{Kuperberg:spiders} G. Kuperberg, \textit{Spiders for rank 2 {Lie} algebras}, Comm. Math. Phys. 180(1) (1996), 109--151.

\bibitem{Maschke:cayley} H. Maschke, \textit{The Representation of Finite Groups}, Amer. J. Math. 16 (1896), 156--194.

\bibitem{Morrison} S. Morrison,
\textit{A Diagrammatic Category for the Representation Theory of $U_q(sl_n)$}, UC Berkeley Ph.D. thesis, arXiv:0704.1503.

\bibitem{mu:alexander} K. Murasugi, \textit{On periodic knots}, Comment. Math. Helv.,
46 (1971), 162--174.

\bibitem{mu:jones} K. Murasugi, \textit{The Jones polynomials of periodic links},
Pacific J. Math., 131 (1988), 319--329.

\bibitem{MOY:Homfly} H. Murakami and T. Ohtsuki and S. Yamada, \textit{HOMFLY polynomial via an invariant of colored plane graphs},
L'Enseignement Mathematique, t., 44 (1998), 325--360.

\bibitem{MT:alternating} W. Menasco and M. Thistlethwaite, \textit{The Classification of Alternating Links}, Ann. Math. 138 (1993), 113--171.

\bibitem{Przytycki:criterion} J. H. Przytycki, \textit{On Murasugi's and Traczyk's criteria for
periodic links}, Math. Ann. 283 (1989), 465--478.

\bibitem{Riley} R. Riley, \textit{An Elliptic Path from Parabolic Representations to Hyperbolic Structures}, In Topology of Low-Dimensional Manifolds, Proceedings, Sussex 1977 (Ed. R. Fenn). New York: Springer-Verlag, pp. 99-133, 1979.

\bibitem{PSTricks}
T.~{Van Zandt}. PSTricks: {PostScript} macros for generic {\TeX}.
Available at {\tt ftp://ftp.princeton.edu/ pub/tvz/}.

\bibitem{Thomassen} C. Thomassen, \textit{Tilings of the torus and the Klein bottle and vertex-transitive graphs on a fixed surface},
Trans. Amer. Math. Soc. 323 (1991), 605--635.

\bibitem{traczyk:period3} P. Traczyk, \textit{A criterion for knots of period $3$}, Topology and its
Appl. 36 (1990), 275--281.

\bibitem{Whitten} W. Whitten, \textit{Symmetries of links}, Trans. Am. Math. Soc. 135 (1969), 213--222.

\bibitem{yokota:skein} Y. Yokota, \textit{The skein polynomial of periodic knots}, Math. Ann.
291(2) (1991), 281--291.

\bibitem{yokota:kauffman} Y. Yokota, \textit{The Kauffman polynomial of periodic knots}, Topology
32(2) (1993), 309--324.

\bibitem{yokota:skeinforn} Y. Yokota, \textit{Skein and quantum
$SU(N)$ invariants of 3-manifolds}, Math. Ann., 307 (1997),
109--138.

\end{thebibliography}
\end{document}